\documentclass[twoside,12pt]{article}
\makeatletter
\def\ifundefined#1{\expandafter\ifx\csname#1\endcsname\relax}
\RequirePackage{theorem}

\AtBeginDocument{
   \makeatletter
   \pagestyle{myheadings}
   \markboth{\hfill\ifundefined{@authorshort} \@author
                   \else \@authorshort 
                   \fi\hfill}
            {\hfill\ifundefined{@titleshort} \@title 
                   \else \@titleshort
                   \fi\hfill}
   \makeatother
 }
     \theorembodyfont{\slshape}
        \newtheorem{thm}{Theorem}[section]

     \theorembodyfont{\upshape}

     \newtheorem{prob}[thm]{Problem}

     \newtheorem{rem}[thm]{\mdseries\scshape Remark}
\ifundefined{proof}

\newenvironment{proof}[1][\proofname]{\par
  \normalfont
  \topsep6\p@\@plus6\p@ \trivlist
  \item[\hskip\labelsep\scshape
    #1{.}]\ignorespaces
}{
  $\qed$\endtrivlist
}
\newcommand{\proofname}{Proof}
\fi

\providecommand{\dedicatory}[1]{}
\providecommand{\keywords}[1]{\begingroup \def \protect {\noexpand \protect \noexpand }\xdef \@thefnmark { }\endgroup \@footnotetext{{\em Keywords and phrases.\/} #1}}
\providecommand{\AMSMSC}[2]{\begingroup \def \protect {\noexpand \protect \noexpand }\xdef \@thefnmark { }\endgroup \@footnotetext{{1991 \it 
Mathematical Subject Classification.\/} Primary: #1; Secondary: #2.}}

\newcommand{\authorshort}[1]{\gdef\@authorshort{#1}}
\newcommand{\titleshort}[1]{\gdef\@titleshort{#1}}
   {\@definecounter{equation}}{\@newctr {equation}[section]}
   
\def\p@enumi{\thethm.}

\def\p@enumi{}
\makeatother

\newcommand{\comment}[1]{}

\usepackage{amsfonts}

\newcommand{\algebra}[1]{\ensuremath{\mathfrak{#1}}}

\newcommand{\Space}[2]{\ensuremath{ {\mathbb{#1}^{#2}} }}

\newcommand{\FSpace}[2]{{\ensuremath{ #1_{#2} }}}

\ifundefined{qed}
    \DeclareMathSymbol{\qed}{0}{AMSa}{"03}
\fi

\newcommand{\scalar}[2]{\langle #1,#2\rangle}

\providecommand{\eqref}[1]{\textup{(\ref{#1})}}

\newcommand{\person}[1]{\textsc{#1}}

\titleshort{A Paley-Wiener Theorem for Nilpotent Lie Groups}
\authorshort{Vladimir V. Kisil}

\begin{document}
\title{A Paley-Wiener Like Theorem\\ for Nilpotent Lie Groups}
\author{Vladimir V. Kisil\\
   Institute of Mathematics\\
			Economics and Mechanics\\
			Odessa State University\\
			ul. Petra Velikogo, 2\\
			Odessa-57, 270057, UKRAINE}
\maketitle
\centerline{\emph{Dedicated to the memory of G.~Polya}}
\begin{abstract}
A version of Paley-Wiener like theorem for connected,
simply connected nilpotent Lie groups is proven.

\end{abstract}
\newpage
\section{Introduction}
In the paper~\cite{Park95} a version of Paley-Wiener theorem for 
two- and three-step nilpotent Lie group was proven. To this end
\person{R.~Park} developed a subtle technique for analysis on
nilpotent Lie groups. It seems that the technique is of separate
interest and could be used to study other problems.

However, it is possible to give a shorter proof of a more general 
theorem based on the completely standard results about nilpotent Lie 
groups. This is the goal of the present paper. We will prove the 
result of \person{Park} for all connected, simply connected 
(exponential) nilpotent Lie groups, but even this is not last level of 
generality---see Remark~\ref{re:general}.
It turns to be that the theorem in such a form is a direct 
consequence of the one-dimensional Paley-Wiener theorem.
 It seems that this
is another example of \emph{inventor paradox}~\cite{Polya62}.

In Section~\ref{se:prelim} we give standard facts about nilpotent
Lie groups, which will be used in Section~\ref{se:theorem} to prove a
version of Paley-Wiener theorem.

\section{Preliminaries}\label{se:prelim}

The following information about nilpotent Lie group could be found in
original papers~\cite{Kirillov62,Kirillov67} or in
monographs~\cite{Kirillov76} and~\cite[Chap.~6]{MTaylor86}. 
It should be noted that the
method of orbits and the induced representations technique are closely
connected for nilpotent Lie groups.

Let $\algebra{g}$ be a nilpotent Lie algebra of the dimension $n$, let
$G=\exp(\algebra{g})$ be the connected, simply connected nilpotent Lie
group corresponding to it. We will identify \algebra{g} and $G$ via
the exponential map and will consider both of them coinciding with 
$\Space{R}{n}$ (as vector spaces and a $C^\infty$-manifolds 
respectively). 
We also will denote by a common letter a representation $\pi$ of $G$ 
and its derived representation of $\algebra{g}$.
All irreducible unitary representations
could be constructed by the inductive 
procedure~\cite{Kirillov62}, \cite[Chap.~6]{MTaylor86}. 
We have

\begin{enumerate}

\item\label{it:beg} The group $G$ is unimodular, the two-sided Haar
measure $d\mu$ on $G$ coincides with the Lebesgue measure 
on $\Space{R}{n}$. Any Lie subgroup $H$ is also nilpotent and
homeomorphic to $\Space{R}{m}$ for some $m\leq n$. The homogeneous
space $X=H\backslash G$ is homeomorphic to $\Space{R}{n-m}$.
Invariant measures on $H$ and $X$ coincide again with
the Lebesgue measures on $\Space{R}{m}$ and $\Space{R}{n-m}$
correspondingly. 

\item\label{it:dual} The unitary dual $\widehat{G}$ could
be parametrized by orbits of co-adjoint representation in
$\algebra{g}'$. 
The support $\widetilde{G}$ of the Plancherel measure $d\nu$ in
$\widehat{G}$ corresponds to orbits of maximal dimensionality 
and parametrized by $\Space{R}{k}\subset \algebra{g}'$, 
where $\dim \algebra{g}=n$ and
$n-k$ is the maximal dimensionality of orbits. 
Moreover~\cite[\S~
7.4]{Kirillov62}, the Plancherel measure is 
equivalent to the 
Lebesgue measure on $\widetilde{G}\cong\Space{R}{k}$:
\begin{equation}\label{eq:plancherel}
d\mu=R(\lambda)\,d\lambda_1 \ldots d\lambda_k.
\end{equation}

\item Every unitary irreducible representation $\pi$ of
$G$ is induced by a one-dimensional representation $\pi_0$ of
a subgroup $H\subset G$~\cite[Theorem~5.1]{Kirillov62}.
The last representation has a form
\begin{equation}\label{eq:repr}
\pi_0(\exp A)=\exp(i\scalar{l}{A}), \qquad A\in \algebra{h} 
\end{equation} 
for some $l\in \algebra{g}'$.

\item\label{it:ind} According to the general scheme of induced 
representation (see~\cite[\S~
13.2]{Kirillov76}, 
\cite[Chap.~5]{MTaylor86}) the representation could be realized as 
follows. Let $\FSpace{L}{2}(G,H, \pi_0)$ be the space
of functions on $G$ with the property
\begin{displaymath}
F(hg)=\pi_0(h)F(g), \qquad h\in H,\ g\in G.
\end{displaymath}
Then $\pi$ is equivalent to the representation on $\FSpace{L}{}(G,H, 
\pi_0)$ by the shift:
\begin{equation}\label{eq:form1}
[\pi(g)F](g_1)= F(g_1 g).
\end{equation}
There is an alternative representation. Let $X=H\backslash G$ be a 
homogeneous space, $s: X \rightarrow G$ be a measured mapping, such  
that $s(Hg)\in Hg$. Define an isomorphism 
$\FSpace{L}{2}(G,H,\pi)\rightarrow\FSpace{L}{2}(X)$ as
\begin{equation}\label{eq:connect}
f(x)=F(s(x)),\qquad  F(hs(x))=\pi_0(h)f(x),\qquad x\in X,\ h\in H.
\end{equation}
Then
\begin{equation}\label{eq:form2}
[\pi(g)f](x)= A(g,x)f(xg),
\end{equation}
where
\begin{displaymath}
A(q,x)=\pi_0(h),\qquad hs(xg)=s(x) g.
\end{displaymath}
\end{enumerate}

\section{A Paley-Wiener Theorem}\label{se:theorem}
\begin{thm}
Let $G=\exp(\algebra{g})$ be a connected, simply connected nilpotent
Lie group. Let $\phi(g)$ be a function from
$\FSpace{L}{\infty}( G)$ with a compact support. Then if 
$\widehat{\phi}(\pi)=0$ on a subset $E$ of $\widehat{G}$
of a positive Plancherel measure then $\phi(g)=0$ 
almost everywhere on $G$.
\end{thm}
\begin{proof}
We start from the change of variables in the formula defines
the non-commutative Fourier transform $\widehat{\phi}(\pi)$ of $\phi$
in the form~\eqref{eq:form1} of induced representations:
\begin{eqnarray*}
[\widehat{\phi}(\pi)F](g_1)&=& \int_G \phi(g) F(g_1 g')\, d\mu(g')\\
  &=& \int_G \phi(g_1^{-1}g) F(g)\, d\mu(g).
\end{eqnarray*}
We rewrite the last line for the representation~\eqref{eq:form2}  
using connection~\eqref{eq:connect}:
\begin{eqnarray}
[\widehat{\phi}(\pi)f](x_1)&=&\int_G \phi([s(x_1)]^{-1} hs(x)) 
    \pi_0(h)f(x)\,d\mu(hs(x))\label{eq:first}\\
  &=&\int_X \int_H \phi([s(x_1)]^{-1} hs(x)) \pi_0(h)f(x)\,dh\,dx 
      \nonumber\\
   &=&\int_X \int_H \phi([s(x_1)]^{-1} hs(x)) \pi_0(h)\,dh\,f(x)\,dx
       \nonumber\\
   &=&\int_X \left( \int_H \phi([s(x_1)]^{-1} hs(x)) 
      \exp(i\scalar{l}{\log h})\,dh\right) f(x)\,dx. \label{eq:last}
\end{eqnarray}
(We substitute in~\eqref{eq:last} the form~\eqref{eq:repr} of 
representation $\pi_0$.) Thus $\widehat{\phi}(\pi)$ 
is an integral operator
$K(l):\FSpace{L}{2}(X)\rightarrow \FSpace{L}{2}(X)$ of the form
\begin{displaymath}
[K(l)f](x_1)=\int_X K(l,x_1,x)f(x)\,dx
\end{displaymath}
with the kernel $K(l,x_1,x)$ defined via the usual Fourier transform 
with respect to variables $h\rightarrow l$:
\begin{displaymath}
K(l,x_1,x)=\int_H \phi([s(x_1)]^{-1} hs(x)) \exp(i\scalar{l}{\log 
h})\,dh.
\end{displaymath}
If the operator $K(l)$ is equal to $0$ for a fixed $l$, then
the kernel $K(l,x_1,x)$ should be equal to zero almost everywhere 
for $(x_1,x)\in X\times X$. But if the last statement fulfils
for all $l\in E$, where $E\subset\Space{R}{k}$ has a non-zero measure,
then $\phi(g)=0$ a.e. by the standard Paley-Wiener theorem.
\end{proof}
\begin{rem}\label{re:general}
To prove the theorem we have used only
properties~\ref{it:beg}--\ref{it:ind}. 
It is well known~\cite{Kirillov62}, \cite[Chap.~6]{MTaylor86} that 
these properties are not too specific and shared also, for example, by 
type I solvable Lie group. 
So it is naturally to assume that our theorem and its proof could
be also transformed to this more general case. 
Reader could found a version of our 
transformations~\eqref{eq:first}--\eqref{eq:last} 
for the group $SL(2,\Space{R}{})$ in~\cite[\S~III.4]{Lang85}.
\end{rem}
We are tempted to conclude the paper by the following
\begin{prob}
Let $G$ be a nilpotent (or solvable) Lie group. Found conditions for 
the function $\phi$, which guarantees that its Fourier transform 
$\widehat{\phi}(\pi)$ is invertible 
almost everywhere in the Plancherel measure on $\widehat{G}$.
\end{prob}

\section{Acknowledgments}
The author was partially supported by the
INTAS grant 93--0322.
During the preparation of the paper in fall of 1995 the author enjoyed
the hospitality of
Universiteit Gent, Vakgroep Wiskundige Analyse, Belgium.
The author is also very grateful to D.S.~Kalyuzhny and the unknown
referee of \textsf{Journal of Functional Analysis} for 
useful comments.

\newcommand{\noopsort}[1]{} \newcommand{\printfirst}[2]{#1}
  \newcommand{\singleletter}[1]{#1} \newcommand{\switchargs}[2]{#2#1}
  \newcommand{\irm}{\textup{I}} \newcommand{\iirm}{\textup{II}}
  \newcommand{\vrm}{\textup{V}}

\end{document}